\documentclass{amsart}
 
\usepackage{amssymb,latexsym,amsmath,amsfonts,hhline,latexsym}
\usepackage{pdfsync}

\newcommand{\cal}{\mathcal}
\usepackage{cite}
\usepackage[matrix,arrow,curve]{xy}
\input xy
\xyoption{all}

\oddsidemargin=\evensidemargin \advance\topmargin by-\headheight

\DeclareMathOperator{\aut}{Aut}

\DeclareMathOperator{\cay}{Cay}

\DeclareMathOperator{\fmult}{OMult}

\DeclareMathOperator{\id}{id}

\DeclareMathOperator{\iso}{Iso}

\DeclareMathOperator{\mult}{Mult}

\DeclareMathOperator{\rad}{rad}
\DeclareMathOperator{\rk}{rk}

\DeclareMathOperator{\sym}{Sym}

\def\mQ{{\mathbb Q}}

\def\mZ{{\mathbb Z}}

\def\frS{{\mathfrak S}}

\def\fS{{\Sigma}}

\def\cA{{\cal A}}

\def\cP{{\cal P}}

\def\cS{{\cal S}}

\def\wh{\widehat}

\def\proof{{\bf Proof}.\ }
\def\bull{\vrule height .9ex width .8ex depth -.1ex }
\def\qaq{\quad\text{and}\quad}
\newcommand{\grp}[1]{\langle {#1}\rangle}

\makeatletter
\renewcommand{\subsection}{\@startsection{subsection}{2}{0mm}{-2mm}{-2mm}{\bf\normalsize}}

\def\sbsnt#1{\subsection{#1}}
\makeatother

\newtheorem{formula}{}[section]
\newtheorem{proposition}[formula]{Proposition}
\newtheorem{definition}[formula]{Definition}
\newtheorem{corollary}[formula]{Corollary}
\newtheorem{remark}[formula]{Remark}
\newtheorem{lemma}[formula]{Lemma}
\newtheorem{theorem}[formula]{Theorem}

\newtheorem{example}[formula]{Example}

\def\thrm{\begin{theorem}}
\def\thrml#1{\begin{theorem}\label{#1}}
\def\ethrm{\end{theorem}}
\def\prpstn{\begin{proposition}}
\def\prpstnl#1{\begin{proposition}\label{#1}}
\def\eprpstn{\end{proposition}}
\def\rmrk{\begin{remark}}
\def\rmrkl#1{\begin{remark}\label{#1}}
\def\ermrk{\end{remark}}
\def\dfntn{\begin{definition}}
\def\dfntnl#1{\begin{definition}\label{#1}}
\def\edfntn{\end{definition}}
\def\nmrt{\begin{enumerate}}
\def\enmrt{\end{enumerate}}
\def\tm#1{\item[{\rm (#1)}]}
\def\qtn{\begin{equation}}
\def\qtnl#1{\begin{equation}\label{#1}}
\def\eqtn{\end{equation}}
\def\lmm{\begin{lemma}}
\def\lmml#1{\begin{lemma}\label{#1}}
\def\elmm{\end{lemma}}
\def\crllr{\begin{corollary}}
\def\crllrl#1{\begin{corollary}\label{#1}}
\def\ecrllr{\end{corollary}}
\def\hpthss{\begin{hypothesis}}
\def\hpthssl#1{\begin{hypothesis}\label{#1}}
\def\ehpthss{\end{hypotxesis}}
\def\xmpl{\begin{example}}
\def\xmpll#1{\begin{example}\label{#1}}
\def\exmpl{\end{example}}
\def\css{\begin{cases}}
\def\ecss{\end{cases}}

\begin{document}
\title[On separability problem for circulant S-rings]{On separability problem for circulant S-rings} 
\author{Sergei Evdokimov}
\address{St. Petersburg Department of Steklov Institute of Mathematics, Russia}
\email{evdokim@pdmi.ras.ru}
\author{Ilya Ponomarenko}
\address{St. Petersburg Department of Steklov Institute of Mathematics, Russia}
\thanks{The second author was partially supported by RFFI Grant 14-01-00156}
\email{inp@pdmi.ras.ru}

\begin{abstract}
A Schur ring (S-ring) over a group $G$ is called separable if every of its similaritities  is induced by isomorphism.
We establish a criterion for an S-ring to be separable in the case when the group $G$ is cyclic. Using
this criterion, we prove that any S-ring over a cyclic $p$-group is separable and that
the class of separable circulant S-rings is closed with respect to duality.
\end{abstract}

\maketitle

\section{Introduction}
A Schur ring or S-ring over a finite group $G$ can be defined as a subring of the group ring $\mZ G$ that is a free 
$\mZ$-module spanned by a partition of $G$ closed under taking inverse and containing $\{1_G\}$ as a class (see 
\cite{EP14} for details). 
 Here a class  of the partition or {\it basic set}, is also treated as the sum of its elements in the group ring.
A {\it Cayley isomorphism} of two S-rings is defined as an isomorphism of the underlying
groups that induces a ring isomorphism of them; a {\it similarity} of two S-rings is defined as a ring
isomorphism that induces a bijection between their basic sets. Any Cayley isomorphism
induces in a natural way a similarity, whereas not every similarity is induced by a Cayley isomorphism.\medskip

The third type of S-ring isomorphisms comes from graph theory. Namely, under {\it isomorphism} of 
S-rings $\cA$ and $\cA'$ over groups $G$  and $G'$ respectively, we mean a bijection $f:G\to G'$ such that
for any basic set $X$ of $\cA$ the image of the Cayley graph $\cay(X,G)$ with respect to $f$, is
the Cayley graph $\cay(X',G')$  where $X'$ is a basic set of $\cA'$.  One can prove that the  bijection $X\mapsto X'$ 
induces by  linearity a similarity $\varphi$ from $\cA$ onto $\cA'$ such that
\qtnl{260613a}
f(Xy)=X^\varphi f(y),
\eqtn
for all basic sets $X$ of $\cA$ and  $y\in G$ where $X^\varphi=X'$.\footnote{In fact, the isomorphism of S-rings
can be defined by means of~\eqref{260613a}, see~\cite{EP14}.} In this case, we say that $\varphi$
is {\it induced} by~$f$. For a fixed similarity~$\varphi$, 
the set of all such isomorphisms $f$ is denoted by $\iso(\cA,\cA',\varphi)$.
Every Cayley isomorphism is also an isomorphism, but the converse statement is not true.\medskip

The {\it separability problem} for S-rings can be formulated as follows: given S-rings $\cA$ and $\cA'$,
and a similarity $\varphi:\cA\to\cA'$ check whether 
$$
\iso(\cA,\cA',\varphi)\ne\varnothing. 
$$
From the computational
complexity point of view this problem is polynomial-time equivalent to the isomorphism problem for Cayley graphs~\cite{EP09a}.
One of the most interesting questions connected with the separability problem is to identify the {\it separable} S-rings~$\cA$,
i.e. those for which the set $\iso(\cA,\cA',\varphi)$ is not empty for all~$\varphi$. In the present paper  we do 
this for circulant  S-rings  (i.e.  when the groups $G$ and $G'$ are cyclic).\medskip

The importance of separable S-rings comes from the observation that any separable S-ring is determined up to isomorphism
by the array of structure constants (with respect to the basis corresponding to the partition of the underlying group). 
However, even for circulant S-rings the situation is quite complicated: there are infinitely many both separable and 
non-separable circulant S-rings (see \cite[Theorem~6.6]{EP01ce}   and \cite[Theorem~1.1]{EP01}). To characterize the
separable circulant S-rings, the theory of circulant S-rings developed by the authors in a series of papers will be
considerably used. Let us briefly discuss its relevant notions.\medskip

Let $\cA$ be an S-ring over a cyclic group $G$. In proving that $\cA$ is separable, without loss of generality we
can assume that $G'=G$, and also that $\cA'=\cA$ \cite{M94}. Next, for technical reasons it is convenient
to restrict oneself to the case when $\cA$ is a {\it quasidense} S-ring, i.e. when it has no rank~$2$ sections of composite
order. This restriction is not essential, because every S-ring over $G$ is separable if and only if so is every quasidense 
S-ring over $G$ (Theorem~\ref{071014a}).\medskip

The proofs in this paper are heavily based on the theory of coset S-rings developed in \cite{EP14}; recall that an S-ring 
over $G$ is a {\it coset} one if  any of its basic sets is a coset of a subgroup 
in~$G$. It was proved there that any circulant coset S-ring is schurian and separable. Our first theorem gives
a criterion for a quasidense S-ring~$\cA$ to be separable, in terms  of  its {\it coset closure} $\cA_0$ which is by definition
the intersection of all coset S-rings containing $\cA$ (it should be noted that in our case $\cA_0$ is itself
a coset S-ring).\medskip

Below we denote by $\Phi(\cA)$ and $\Phi_\infty(\cA)$  the group of all similarities of $\cA$
and the group of all $\varphi\in\Phi(\cA)$ with $\iso(\cA,\varphi)\ne\varnothing$ where
$\iso(\cA,\varphi)=\iso(\cA,\cA,\varphi)$.

\thrml{280414a}
Let $\cA$ be a quasidence circulant S-ring and $\cA_0$ its coset closure. Then
$$
\Phi_\infty(\cA)=\Phi(\cA_0)^\cA
$$
where $\Phi(\cA_0)^\cA$ is the subgroup of $\Phi(\cA)$ induced by the action of~$\Phi(\cA_0)$ on~$\cA$.
In particular, $\cA$  is separable if and only if $\Phi(\cA)=\Phi(\cA_0)^\cA$.
\ethrm

To apply the separability criterion given in Theorem~\ref{280414a} one has to have on hand the S-ring $\cA_0$ and the
groups $\Phi(\cA)$ and  $\Phi(\cA_0)$. However, they cannot in general  be easily found. That is why we would like to
simplify this criterion by using the observation that every similarity $\varphi\in\Phi(\cA)$ is uniquely determined
by its restrictions to principal $\cA$-sections, and hence to sections belonging to  the class $\frS_0(\cA)$ of all $\cA$-sections 
of~$G$ that are projectively equivalent to subsections of principal $\cA$-sections. This enables us to replace $\varphi$
by the family of these restrictions, and eventually to get an explicit  form of the criterion that is  in principle can be
reduced to the compatibility of a modular linear system as in~\cite{EP14}.\medskip

Let $\cA$ be a circulant S-ring. Suppose that for  every section $S\in\frS_0(\cA)$ we are given an automorphism
$\sigma_S$ of the group~$S$. Then the family
$$
\fS=\{\sigma_S\}_{S\in\frS_0(\cA)}
$$
is called a {\it $\cA$-multiplier} if for any sections $S,T\in\frS_0(\cA)$
such that $T$ is projectively equivalent to a subsection of~$S$, the automorphisms $\sigma_T\in\aut(T)$ and
$\sigma_S\in\aut(S)$ are induced by raising to the same power\footnote{We recall that any
automorphism of a finite cyclic group is induced by raising to a power coprime to the order of this group.}
(see Definition~\ref{200313b}). The set of all of them forms a subgroup $\mult(\cA)$ of the direct product 
$\prod_{S\in\frS_0(\cA)}\aut(S)$.\footnote{The multipliers of $\cA$ defined in~\cite{EP14} form a subgroup of
the group $\mult(\cA)$; it seems reasonable to call them {\it inner multipliers.}}
\medskip

For an $\cA$-section $S$ denote by $\aut_\cA(S)$ 
the group of Cayley automorphisms of the restriction $\cA_S$ of the S-ring~$\cA$ to~$S$.
Then every  $\cA$-multiplier $\fS$ produces a family 
$$
\fS'=\{C_S\}_{S\in\frS_0(\cA)}
$$ 
where  $C_S=\aut_\cA(S)\sigma_S$. For this family the compatibility condition similar to that as for a $\cA$-multiplier,
is obviously satisfied. Therefore 
in the sense of Definition~\ref{200313bx}, the family $\fS'$ is an {\it outer $\cA$-multiplier}. It follows that the mapping
\qtnl{260814b}
\theta:\mult(\cA)\to \fmult(\cA),\quad \fS\mapsto\fS'
\eqtn
is a group homomorphism where $\fmult(\cA)$ is the group of all outer $\cA$-multipliers. In these terms the explicit
form of the criterion looks as follows.

\thrml{170113a}
A quasidense circulant S-ring $\cA$ is separable if and only if  the homomorphism $\theta$ is surjective.
\ethrm

From the Burnside theorem on permutation groups of prime degree it follows that every S-ring over a group  of
prime order $p$, is normal or of rank~2. In both cases it is separable. Therefore all S-rings over such a group
are separable. Using Theorem~\ref{170113a} we show that any cyclic $p$-group has this property.


\thrml{080914a}
Any  S-ring over a cyclic $p$-group is separable.
\ethrm

There is a deep relationship between the schurity and separability problems (see e.g. \cite{EP09a}). 
For circulant S-rings the first problem was solved in~\cite{EP14}.  The schurity criterion obtained there
implies that the class of schurian circulant S-rings is closed with respect to duality. The analog of this result for
separable circulant S-rings is also true.

\thrml{170113b}
A circulant S-ring is separable if and only if so is the dual  S-ring.
\ethrm

We complete the introduction by remarking that this paper is closely related to paper~\cite{EP14}
by techniques, results and spirit. Therefore here we follow the notation and terminology of that paper. When referring to it,
we keep only the number of the statement, preceding it by the letter~A (e.g. instead of \cite[Theorem~4.2]{EP14}
we write Theorem~A4.2). 

\section{Auxiliary statements on S-rings}

\sbsnt{Similarities and isomorphisms.} The basics on isomorphisms and similarities of S-rings can be found in
Section~A3. Below we are mainly interested in the question how to extend them.

\lmml{210714d}
Let $\cA_i$ be an S-ring over a group $G_i$ and let $\cA'_i$ be the minimal S-ring over $G_i$ that contains $\cA_i$ and
an element $\xi_i\in\mZ G$, $i=1,2$. Then given a similarity $\varphi\in\Phi(\cA_1,\cA_2)$ there is at most one
similarity $\varphi'\in\Phi(\cA'_1,\cA'_2)$ extending $\varphi$ and such that $\varphi'(\xi_1)=\xi_2$.
\elmm
\proof Every element of $\cA'_i$ can be obtained from the elements of $\cA_i$ and $\xi_i$ by
taking $\mQ$-linear combinations and products (ordinary and entrywise), $i=1,2$. Since the $\mQ$-linear extension
$\psi$ of $\varphi'$ must preserve all these operations, it is uniquely determined (if exists).
This implies the required statement because $\psi$ and $\varphi'$ coincide on $\cA'_1$.\bull\medskip

It follows from equality (A12) that given two isomorphic S-rings one can choose their isomorphism to be {\it normalized},
i.e. preserving the neutral elements (any bijection between two groups that has this property, is also called normalized).
The normalized isomorphisms admit a more natural description than in~\eqref{260613a} (see Section~A3.2).

\lmml{241014a}
Let $\cA_1$ and $\cA_2$ be S-rings over groups $G_1$ and $G_2$ respectively, and let $f:G_1\to G_2$ be a 
normalized bijection. Suppose that  $f$ preserves the basic sets. Then $f\in\iso(\cA_1,\cA_2)$
if and only if $f(Xy)=f(X)f(y)$ for all $X\in\cS(\cA_1)$ and $y\in G_1$.\bull
\elmm

Any isomorphism  $f$ of two association schemes is also isomorphism from the extension of the first scheme by a relation $R$
onto the extension of the second scheme by the relation $R^f$. A similar statement does not hold in general for 
Cayley schemes, and hence for S-rings.  However, it becomes true under  some extra condition.

\lmml{140514a}
Let $\cA_i$ be an S-ring over a group $G_i$ and $X_i\subset G_i$, $i=1,2$.  Suppose that $f\in\iso(\cA_1,\cA_2)$ is 
a normalized isomorphism such that  $f(X_1)=X_2$ and
$$
f(X_1y)=X_2f(y)
$$
for all $y\in G_1$. Then
$f\in\iso(\cA'_1,\cA'_2)$ where $\cA'_i$ is the minimal S-ring over $G_i$ that contains $\cA_i$ and $X_i$.
\elmm
\proof Let us extend $f$ to a linear isomorphism from $\mQ G_1$ to $\mQ G_2$; it will be denoted by the same
letter. Set $M$ to be the set of all elements $\xi\in\mQ G_1$ such that $f(\xi x)=f(\xi)f(x)$ for all $x\in G_1$. Then 
$\cA_1\subset M$ by  Lemma~\ref{241014a}, and $\sum_{x\in X_1}x\in M$ by lemma hypothesis. Besides, obviously 
$M$ is a $\mQ$-module such that $M\circ M\subset M$. Furthermore, if $M$ contains both 
$\xi$ and $\eta=\sum_{x\in G_1}c_xx$, then
$$
f(\xi\eta y)=\sum_{x\in G_1}c_xf(\xi xy)=\sum_{x\in G_1}c_xf(\xi)f(xy)=f(\xi)\sum_{x\in G_1}c_xf(xy)=
$$
$$
f(\xi)f(\eta y)=f(\xi)f(\eta)f(y)=(\sum_{x\in G_1}c_xf(\xi)f(x))f(y)=
$$
$$
(\sum_{x\in G_1}c_xf(\xi x))f(y)=f(\xi\eta)f(y)
$$
for all $y\in G_1$. Thus $M\, M\subset M$. Since every element of $\cA'_1$ can be obtained from $\sum_{x\in X_1}x$ and
the elements of $\cA_1$   by taking $\mQ$-linear combinations and products (ordinary and entrywise),
we conclude that $\cA'_1\subset M$. On the other hand, $f(\cA'_1)$ is the $\mZ$-module associated with
the partition $f(\cS(\cA_1))$ of the group $G_2$ where $\cS(\cA_1)$ is the partition of $G_1$ into the basic sets of~$\cA_1$.
Besides, $f(\cA'_1)$ is closed with respect to multiplication because $\cA'_1\subset M$. So this module 
is an S-ring over $G_2$. However, it contains both $\cA_2=f(\cA_1)$ and $\sum_{x\in X_2}x=f(\sum_{x\in X_1}x)$.
Thus by the minimality of $\cA'_2$ we have $f(\cA'_1)\supset\cA'_2$. Similarly, interchanging $\cA_1$ and $\cA_2$
and taking $f^{-1}$ instead of $f$, one can prove that $f^{-1}(\cA'_2)\supset\cA'_1$.  Therefore
$f(\cA'_1)=\cA'_2$, and hence $f$ is a linear isomorphism from $\cA'_1$ onto $\cA'_2$.  Thus $f\in\iso(\cA'_1,\cA'_2)$,
as required.\bull\medskip

Any section $U/L$ of a group $G$ naturally acts on the right $L$-cosets in $G$, namely $Lx\cdot Ly:=Lxy$
for $x\in U$ and $y\in G$. If $U/L$ is a section of an S-ring $\cA$ over $G$, then a normalized isomorphism $f\in\iso(\cA)$
takes this section to another $\cA$-section, but in general does not preserve the above action. Condition~\eqref{140115a}
in the lemma below, means that $f$ does preserve it.

\crllrl{190514e}
Let $\cA$ be an S-ring over a group $G$, $S=U/L$ an $\cA$-section and let $f\in\iso(\cA)$ be a normalized isomorphism.
 Suppose that 
\qtnl{140115a}
f(Lxy)=f(Lx)f(Ly)
\eqtn
for all $x\in U$ and $y\in G$. Then $f\in\iso(\cA')$ where
$\cA'$ is the minimal S-ring over $G$ that contains $\cA$ and for which $(\cA')_S=\mZ S$.
\ecrllr
\proof  The S-ring $\cA'$ is obtained from $\cA$ by subsequent extensions of $\cA$ by means of $L$-cosets in~$U$.
Therefore it suffices to verify that given $x\in U$, $f$ is an isomorphism of  the minimal S-ring over $G$ that contains
$\cA$ and $Lx$.  However, this follows from Lemma~\ref{140514a} for $\cA_1=\cA_2=\cA$ and $X_1=Lx$,
because by the hypothesis on $f$ we have
$$
f(Lxy)=f(Lx)f(Ly)=f(Lx)f(L)f(y)=f(Lx)f(y),
$$
here we used the fact that $f(Ly)$ equals the $f(L)$-coset containing $f(y)$.\bull

\sbsnt{Projective equivalence.}
Let $S=U/L$ and $S'=U'/L'$  be sections of a group $G$. Suppose that $S'$ is a multiple of $S$, i.e.
that $UL'=U'$ and $U\cap L'=L$.
Set $f_{S,S'}:xL\mapsto xL'$ to be the canonical projective isomorphism from $S$ to $S'$ defined in
\cite[Section~3.1]{EP12}; we also set $f_{S,S'}:=(f_{S',S})^{-1}$ when $S$ is a multiple of $S'$. Then
the following statement obviously holds.

\lmml{180614b}
Suppose that $S$ is a multiple of both sections $S'$ and $S''$one of which is a multiple
of the other. Then
$$
f_{S^{},S'}f_{S',S''}=f_{S^{},S''}.\bull
$$
\elmm

Let now the sections $S$ and $S'$ be projectively equivalent. Then there exist sections $S=S_1,S_2,\ldots,S_{k-1},S_k=S'$
such that for all $i=1,\ldots,k-1$ one of the sections $S_i$, $S_{i+1}$ is a multiple of the other. The mapping $f:S\to S'$
defined by
\qtnl{180614c}
f=f_{S_1,S_2}\cdots f_{S_{k-1},S_k},
\eqtn
is called a projective isomorphism from $S$ onto $S'$.

\thrml{180614a}
A projective isomorphism between two projectively equivalent sections of a cyclic group is uniquely determined.
\ethrm
\proof
Let $S$ and $S'$ be projectively equivalent sections and let $f:S\to S'$ be the  projective
isomorphism~\eqref{180614c}. Denote by $S_0$ the largest  section  in
the class of projectively equivalent sections that contains $S_1$ (apply Theorem~A5.4 for the group ring). Then $S_0$ is a multiple
of $S_i$ for all~$i$. By Lemma~\ref{180614b} this implies that $f_{S_0,S_i}f_{S_i,S_{i+1}}=f_{S_0,S_{i+1}}$
for all $i=1,\ldots,k-1$. Thus from~\eqref{180614c} it follows that
$$
f=f_{S_1,S_0}f_{S_0,S_k}.
$$
This proves the required statement.\bull

\crllrl{270814a}
Let $S$ and $T$ be projectively equivalent sections of a cyclic group~$G$. Then any section $S'\preceq S$
is projectively equivalent to the section $T'=(S')^{f_{S,T}}$ where $f_{S,T}$ is the projective isomorphism
from $S$ onto $T$.\footnote{We recall that notation $S'\preceq S$
means that $S'$ is a subsection of $S$.}  Moreover, the restriction of $f_{S,T}$ to $S'$ equals $f_{S',T'}$.
\ecrllr
\proof The statement immediately follows from Theorem~\ref{180614a} when $T$ is a multiple of $S$.  In a general
case,  we are done by induction.\bull\medskip

Let $\cA$ be an S-ring over a group $G$.  Under the projective equivalence "$\sim$" on the set $\frS(\cA)$ we mean 
the transitive  closure of the relation "to be a multiple" on this set. Clearly, any two projectively equivalent $\cA$-sections
are projectively equivalent as sections of~$G$. The converse statement is not true in general, but holds
when the group $G$ is cyclic.

\lmml{170614a}
Let $\cA$ be a circulant S-ring, $S$ and $T$ projectively equivalent $\cA$-sections and $f_{S,T}$
the projective isomorphism from $S$ to $T$.  Then
\nmrt
\tm{1}  if $\varphi$ is a similarity of $\cA$,  then $f_{S,T}$  takes $\varphi_S$ to $\varphi_T$,
\tm{2}  if $f$ is a normalized isomorphism of~$\cA$, then $f_{S,T}$  takes $f^S$ to $f^T$.
\enmrt
\elmm
\proof By the definition of projective isomorphism it suffices to verify   statement~(1) in the
case when $T=U_2/L_2$ is a multiple of $S=U_1/L_1$. Without loss of generality we can assume
that $U_2=G$ and $L_1=1$. By \cite[Theorem~3.2]{EP12} the projective isomorphism $f_{S,T}$  is a Cayley isomorphism 
from $\cA_S$ onto $\cA_T$. Denote by $\psi$ the induced similarity. Then
obviously $X^\psi=\pi(X)$ for any $X\in\cS(\cA_{U_1})$ where $\pi$ is the quotient epimorphism
from $U_2$ onto~$T$. Therefore
$$
(X^{\varphi_S})^\psi=\pi(X^{\varphi_S})=
\pi(X^\varphi)=
\pi(X)^{\varphi_T}=
(X^\psi)^{\varphi_T}
$$
and we are done. Statement~(2) immediately follows from \cite[Lemma~3.1]{EP12} with $\Gamma=\iso(\cA)$.\bull

\sbsnt{Duality.}
Let $G$ be an abelian group and $\wh G$ its dual group, i.e. the group of complex characters of $G$. 
Given $\sigma\in\aut(G)$ denote by $\wh\sigma$ the automorphism of $\wh G$ such that 
$\chi^{\wh\sigma}(g)=\chi(g^\sigma)$
for all $g\in G$ and $\chi\in\wh G$. Clearly, the mapping $\sigma\mapsto\wh\sigma$ induces an isomorphism
from $\aut(G)$ onto $\aut(\wh G)$.

\lmml{010914e}
Let $G$ be a cyclic group. Then given sections $S,T\in\frS(G)$ such that $S\preceq T$ we have
$\wh{\sigma^S}=\wh\sigma^{\wh S}$ for all $\sigma\in\aut(T)$ where $\wh S$ is the section of $\wh G$ dual to $S$. 
\elmm
\proof We observe that $\wh S\preceq\wh T$ by statement~(3) of Lemma A4.4.  Since $G$ and $\wh G$ are
cyclic groups, there are restriction epimorphisms $\pi:\aut(T)\to\aut(S)$ and $\wh\pi:\aut(\wh T)\to\aut(\wh S)$, and
the diagram
\begin{figure}[h]
$\xymatrix@R=10pt@C=20pt@M=0pt@L=5pt{
\aut(T)\ \ar@{->}[r]^\pi\ar@{->}[d] & \ \aut(S) \ar@{->}[d]\\
\aut(\wh T)\  \ar@{->}[r]^{\wh\pi} & \ \aut(\wh S)\\
}$
\end{figure}

\noindent is commutative where the vertical arrows mean the isomorphisms $\sigma\mapsto\wh\sigma$.  Therefore  for all $\sigma\in\aut(T)$, $g\in S$ and $\chi\in\wh S$ we have
$$
\chi(g^{\pi(\sigma)})=\wh{\pi(\sigma)}(\chi(g))=\wh\pi(\wh\sigma)(\chi(g))=\chi^{\wh\pi(\wh\sigma)}(g).
$$
It follows that
$$
\chi^{\wh{\sigma^S}}(g)=\chi(g^{\sigma^S})=\chi(g^{\pi(\sigma)})=\chi^{\wh\pi(\wh\sigma)}(g)
=\chi^{\wh\sigma^{\wh S}}(g),
$$
as required.\bull

\lmml{010914a}
Let $S$ and $T$ be projectively equivalent sections of a cyclic group~$G$, and $f:=f_{S,T}$ and $\wh f:=f_{\wh S,\wh T}$
projective isomorphisms. Then $\wh{\sigma^f}=\wh{\sigma}^{\wh f}$ for all $\sigma\in\aut(S)$.
\elmm
\proof We have to verify that the following diagram is commutative
\begin{figure}[h]
$\xymatrix@R=10pt@C=30pt@M=0pt@L=5pt{
\aut(S)\ \ar@{->}[r]^{\sigma\mapsto\sigma^f}\ar@{->}[d] & \ \aut(T) \ar@{->}[d]\\
\aut(\wh S)\  \ar@{->}[r]^{\wh\sigma\mapsto\wh\sigma^{\wh f}} & \ \aut(\wh T)\\
}$
\end{figure}

\noindent However, this is true because all  isomorphisms of this diagram take an automorphism induced by raising to a 
certain power to the automorphism induced by raising to the same power.\bull

\section{Outer multipliers}\label{200214c}

Let $\cA$ be an S-ring over a cyclic group $G$. Suppose that for every section $S\in\frS_0(\cA)$ we are given a coset
$C_S\subset\aut(S)$  of the group  $\aut_\cA(S)=\aut(\cA_S)\cap\aut(S)$.

\dfntnl{200313bx}
The family $\fS=\{C_S\}_{S\in\frS_0(\cA)}$ is called an outer $\cA$-multiplier if the following two conditions are satisfied for all
sections $S_1,S_2\in\frS_0(\cA)$:
\nmrt
\tm{M1} if $S_1\succeq S_2$, then $(C_{S_1})^{S_2}=C_{S_2}$,
\tm{M2} if $S_1\sim S_2$, then $(C_{S_1})^{f_{S_1,S_2}}=C_{S_2}$
\enmrt
where $f_{S_1,S_2}$ is the projective isomorphism from $S_1$ onto $S_2$.
\edfntn

Clearly, the product of two $\aut_\cA(S)$-cosets is also $\aut_\cA(S)$-coset. This enables
us to define the componentwise product of any two outer $\cA$-multipliers. Since for the corresponding family,
conditions (M1) and (M2) are obviously satisfied,  the set of all outer $\cA$-multipliers  forms
a group; we denote it by $\fmult(\cA)$.\medskip

Let $\fS\in\fmult(\cA)$ and $T\in\frS(\cA)$. Then obviously the class $\frS:=\frS_0(\cA_T)$ is contained in $\frS_0(\cA)$,
and conditions (M1) and (M2) are  satisfied for the family $\fS^T=\{C_S\}_{S\in\frS}$.
This proves the following statement.

\lmml{261113b}
The family $\fS^T$ is an outer $\cA_T$-multiplier.\bull
\elmm

Let $\cA$ be a quasidense circulant S-ring and  let $\varphi\in\Phi(\cA)$ be a similarity. Then
for any section $S\in\frS_0(\cA)$ the similarity $\varphi_S$ is induced by an automorphism $\sigma_S\in\aut(S)$.
Indeed, if $S$ is a principal section, then the S-ring $\cA_S$ has trivial radical, and this follows from
Corollary~6.4 and statement~(1) of Theorem~6.6 of \cite{EP01ce}; if $S$ is a subsection of a principal section, then
$\sigma_S$ is obtained by  restriction; if $S$ is projectively equivalent to a subsection of a principal section,
then $\sigma_S$ is obtained by means of projective isomorphism (Theorem A4.2). Set
$$
\fS(\varphi)=\{C_S(\varphi)\}_{S\in\frS_0(\cA)}.
$$
where $C_S(\varphi)=\aut_\cA(S)\sigma_S$.

\lmml{281113d}
The family $\fS(\varphi)$ is an outer $\cA$-multiplier. Moreover, if $T\in\frS(\cA)$, then
$\fS(\varphi)^T=\fS(\varphi_T)$.
\elmm
\proof To prove the first statement, let $S_1,S_2\in\frS_0(\cA)$.  Then the similarity $\varphi_{S_i}$ is induced by
an automorphism $\sigma_{S_i}$ for $i=1,2$. Therefore if $S_1\succeq S_2$, then  $\varphi_{S_2}$ is induced by
$(\sigma_{S_1})^{S_2}$.  It
follows that
$$
(\sigma_{S_2})^{-1}(\sigma_{S_1})^{S_2}\in\aut_\cA(S_2).
$$
Thus $C_{S_2}(\varphi)=C_{S_1}(\varphi)^{S_2}$ and condition~(M1) is satisfied. Since condition (M2) follows from
statement~(1) of Lemma~\ref{170614a}, $\fS(\varphi)$ is an outer $\cA$-multiplier  as required.  The second statement of the lemma is obvious.\bull\medskip

The following theorem is the main result of this subsection.

\thrml{071212d}
Let $\cA$ be an quasidense circulant S-ring. Then the mapping
\qtnl{221113w}
\Phi(\cA)\to \fmult(\cA),\quad \varphi\mapsto\fS(\varphi)
\eqtn
is a group isomorphism.
\ethrm
\proof
Mapping~\eqref{221113w} is obviously a group homomorphism. To prove its injectivity
suppose that $\fS(\varphi)=\fS(\psi)$ for some similarities $\varphi,\psi\in\Phi(\cA)$. Then obviously
$\varphi_S=\psi_S$ for all sections $S\in\frS_0(\cA)$. Therefore $\varphi$ and $\psi$ are equal on
all principal $\cA$-sections. Thus $\varphi=\psi$ by Lemma~A3.1.\medskip

Let us prove the surjectivity of homomorphism~\eqref{221113w} by induction on the size of the group~$G$
underlying~$\cA$. Let $\fS=\{C_S\}_{S\in\frS_0(\cA)}$ be an outer $\cA$-multiplier. First, suppose that $\rad(\cA)=1$.
Then $G\in\frS_0(\cA)$, the S-ring $\cA$ is cyclotomic, and any $\sigma\in C_G$ is a Cayley isomorphism of $\cA$
(\cite[Theorem~2.4]{EKP13}). Clearly, the similarity induced by $\sigma$, does not depend on the choice of it.
Denote this similarity by  $\varphi$. Then 
$\varphi_S$ is induced by $\sigma^S$ for all $S\in\frS_0(\cA)$. Now
condition~(M1) implies that
$$
C_S(\varphi)=\aut_\cA(S)\sigma^S=(\aut_\cA(G)\sigma)^S=(C_G)^S=C_S.
$$
Thus $\fS(\varphi)=\fS$.\medskip

Now, assume that $\rad(\cA)\ne 1$. Then by Theorem~A5.2 the S-ring $\cA$ is a proper  $U/L$-wreath product.
By the inductive hypothesis applied to the S-ring $\cA_U$ and its outer multiplier $\fS^U$, as well as
to  the S-ring $\cA_{G/L}$ and its outer multiplier $\fS^{G/L}$, there exist similarities
$\varphi_1\in\Phi(\cA_U)$ and $\varphi_2\in\Phi(\cA_{G/L})$ such that
\qtnl{291113b}
\fS^U=\fS(\varphi_1)\qaq\fS^{G/L}=\fS(\varphi_2).
\eqtn
By the first equality in~\eqref{291113b} and Lemma~\ref{281113d} applied to $\cA=\cA_U$, $\varphi=\varphi_1$ 
and $T=U/L$, we have
$$
\fS((\varphi_1)_{U/L})=\fS(\varphi_1)^{U/L}=(\fS^U)^{U/L}=\fS^{U/L}.
$$
Similarly, by the second equality in~\eqref{291113b} and Lemma~\ref{281113d} applied to $\cA=\cA_{G/L}$,
$\varphi=\varphi_2$ and $T=U/L$ we have
$$
\fS((\varphi_2)_{U/L})=\fS(\varphi_2)^{U/L}=(\fS^{G/L})^{U/L}=\fS^{U/L}.
$$
Thus $\fS((\varphi_1)_{U/L})=\fS((\varphi_2)_{U/L})$ and by the injectivity statement we have
$$
(\varphi_1)_{U/L}=(\varphi_2)_{U/L}.
$$
Thus by  statement~(2) of Theorem~A3.3 there exists
a unique similarity $\varphi\in\Phi(\cA)$ such that
\qtnl{110414d}
\varphi_U=\varphi_1\qaq\varphi_{G/L}=\varphi_2.
\eqtn

To complete the proof let us verify that $\fS(\varphi)=\fS$. Let $S\in\frS_0(\cA)$. Suppose first that
$S$ is a principal section. Then it is of trivial radical. Therefore $S$ is either an $\cA_U$- or $\cA_{G/L}$-section.
Since $\fS(\varphi)^U=\fS^U$ and $\fS(\varphi)^{G/L}=\fS^{G/L}$ (see~\eqref{291113b} and~\eqref{110414d}),
we conclude that
\qtnl{260814a}
C_S(\varphi)=C_S.
\eqtn
 Let now $S$ be a subsection of a principal section $T$. Then by above and condition (M1) we obtain that
$$
C_S(\varphi)=C_T(\varphi)^S=(C_T)^S=C_S,
$$
and~\eqref{260814a} holds. Finally,  let $S$ be a projectively equivalent to a subprincipal $\cA$-section~$T$. By what 
we just proved, statement~(1) of Lemma~\ref{170614a} and condition (M2),  we obtain that
$$
C_S(\varphi)=C_T(\varphi)^{f_{T,S}}=(C_T)^{f_{T,S}}=C_S
$$
as required.\bull

\section{Reduction to quasidence S-rings}

The main result in this section is the following theorem reducing the separability problem for circulant
S-rings to the quasidense case.

\thrml{071014a}
Let $G$ be a cyclic group. Then every S-ring over $G$ is separable if and only if so is every quasidense S-ring
over $G$.
\ethrm

Theorem~\ref{071014a} will be deduced from Theorem~\ref{270514a} for the proof of which we
recall,  as in~\cite{EP14}, some facts  from paper~\cite{EP12}. Let $\cA$ be an S-ring over a cyclic group $G$. 
A class  $C$ of projectively equivalent $\cA$-sections is called {\it singular} if its  rank is~$2$, its order is greater than~$2$ and it contains two
sections $L_1/L_0$ and $U_1/U_0$ such that the second is a multiple of the first and
the following two conditions are satisfied:
\nmrt
\tm{S1} $\cA$ is both the $U_0/L_0$- and $U_1/L_1$-wreath product,
\tm{S2} $\cA_{U_1/L_0}=\cA_{L_1/L_0}\otimes\cA_{U_0/L_0}$.
\enmrt
By \cite[Lemma~6.2]{EP12}  the above two sections are necessarily  the smallest and largest $\cA$-sections of $C$.
Moreover, from Theorem~4.6 ibid., it follows
that any rank~$2$ class of  composite order belonging to $\cP(\cA)$ is singular; in particular,  the S-ring $\cA$ is
quasidense if and only if no class in $\cP(\cA)$ is singular.  It is worth mentioning that once $\cA$ has a singular class, 
it admits a lot of non-trivial automorphisms; more exactly, the following statement is implied by Lemma~4.3 of 
paper~\cite{EP03}.

\lmml{041114a}
Let $\cA$ be  an S-ring over a cyclic group $G$ and let $U/L$ be the largest section in a singular class of $\cA$.
Then
$$
\aut(\cA)^{G/L}\ge \prod_{X\in G/U}\sym(X/L).
$$
\elmm

For an $\cA$-section $S$ we define the {\it $S$-extension} of~$\cA$ to be the  smallest
S-ring $\cA'\ge\cA$ such that $\cA'_S=\mZ S$. From Theorem~A4.2 it follows that $\cA'$
does not depend on the choice of $S$ in the class  of $\cA$-sections projectively equivalent to~$S$.
The following statement follows from Lemma A13.1 and statements (E1), (E2) inside it.

\lmml{101213a}
Let $\cA$ be a circulant S-ring, $C\in\cP(\cA)$ a singular class, $S=L_1/L_0$ the smallest section in $C$
and  $\cA'$ the $S$-extension of $\cA$. Then
\nmrt
\tm{1} $\rk(\cA')>\rk(\cA)$,
\tm{2} $\cA'$ is both the $U_0/L_0$- and $U_1/L_1$-wreath product,
\tm{3} $\cA'_{U_0}=\cA^{}_{U_0}$, $\cA'_{G\setminus L_1}=\cA^{}_{G\setminus L_1}$ and
$\cA'_{U_1/L_0}=\mZ S\otimes\cA^{}_{U_0/L_0}$.
\enmrt
\elmm

Now we are ready to prove the basic statement on extending similarities to $S$-extensions where
$S$ belongs to a singular class.

\thrml{270514a}
Let $\cA$ be a circulant S-ring, $C\in\cP(\cA)$ a singular class, $S\in C$ and  $\cA'$ the $S$-extension of $\cA$.
Then given similarities $\varphi\in\Phi(\cA)$ and $\psi\in\Phi(\mZ S)$ there exists a unique
similarity $\varphi'\in\cA'$ extending $\varphi$ and such that $\varphi'_S=\psi$.
Moreover, if $\varphi\in\Phi_\infty(\cA)$, then $\varphi'\in\Phi_\infty(\cA')$.
\ethrm
\proof Without loss of generality we can assume that $S=L_1/L_0$ is the smallest section in $C$.
To prove the first statement let $\varphi\in\Phi(\cA)$ and $\psi\in\Phi(\mZ S)$.  Then from
statement~(3) of Lemma~\ref{101213a} it follows that
$$
\varphi_T\in\Phi(\cA'_T),\quad T\in\{U_0,U_0/L_0,G/L_1\}.
$$
Furthermore, since $\rk(\cA_S)=2$, the similarity  $\psi$ extends $\varphi_S$. However,
$\cA'_S=\mZ S$. Thus, again by statement~(3) of Lemma~\ref{101213a}, the mapping $\psi\otimes\varphi_{U_0/L_0}$
is a similarity of $\cA'_{U_1/L_0}$ that extends $\varphi_{U_1/L_0}$. \medskip

By statement~(2) of Lemma~\ref{101213a} the S-ring $\cA'_{U_1}$ is the $U_0/L_0$-wreath product.
Since the restrictions of $\varphi_{U_0}$ and $\psi\otimes\varphi_{U_0/L_0}$ to $U_0/L_0$ coincide,
we find by statement~(2) of Theorem~A3.3 a similarity $\psi'\in\Phi(\cA_{U_1})$ such that
$$
\psi'_{U_0}=\varphi_{U_0}\qaq \psi'_{U_1/L_0}=\psi\otimes\varphi_{U_0/L_0}.
$$
Clearly, $\psi'$ extends $\varphi_{U_1}$ and $\psi'_{U_1/L_1}=\varphi_{U_1/L_1}$.
Furthermore, $\cA'$ is the $U_1/L_1$-wreath product by statement~(2) of Lemma~\ref{101213a}.
Therefore in the same way as above, we find a similarity $\varphi'\in\Phi(\cA')$ such that $\varphi'_{U_1}=\psi'$ and
$\varphi'_{G/L_1}=\varphi_{G/L_1}$.  Since obviously $\varphi'$ extends $\varphi$, we are done.
The uniqueness is clear by Lemma~\ref{210714d} and the definition of $S$-extension.\medskip

To prove the second statement let $\varphi\in\Phi_\infty(\cA)$. Then there exists a normalized isomorphism
$f\in\iso(\cA,\varphi)$. Without loss of generality we can assume that $S=:U/L$ is the largest section in the class~$C$.
For each $Y\in G/U$ fix a permutation $h_Y\in G_{Y\to U}$ with $h_U=\id$.
Besides, by Lemma~A3.2  there exists an automorphism $\sigma\in\aut(S)$ belonging to $\iso(\mZ S,\psi)$.
Then using Lemma~\ref{041114a} one can  choose the isomorphism $f$ so that
$$
f^{Y/L}=h^{Y/L}\cdot\sigma\cdot (h')^{U/L},\qquad Y\in G/U
$$
where $h=h_Y$ and $h'=(h_{Y^f})^{-1}$.  However, the bijections $h$ and $h'$ are induced by multiplications by
elements of~$G$, say $a$ and $a'$, respectively.  Thus the hypothesis of Corollary~\ref{190514e} is satisfied:
indeed, if $x\in U$ and $y\in Y$, then $ay\in U$ and $(Lx)^\sigma=f(Lx)$ (because $h_U=\id$), and hence
$$
f(Lxy)=f^{Y/L}(Lxy)=(Laxy)^\sigma a'=(Lx)^\sigma(Lay)^\sigma a'=
$$
$$
f(Lx)f^{Y/L}(Ly)=f(Lx)f(Ly).
$$
Thus $f\in\iso(\cA')$. By the uniqueness of $\varphi'$ we conclude that $f\in\iso(\cA',\varphi')$
as required.\bull\medskip

The first of two statements below immediately follows from Theorem~\ref{270514a} and is used in the proof
of Theorem~\ref{071014a}; the second one will be used in subsequent sections.

\crllrl{280514f}
In the notation of Theorem~\ref{270514a} the S-ring~$\cA$ is separable if and only if so is $\cA'$.\bull
\ecrllr

\crllrl{280514a}
Let $\cA$ be a circulant S-ring. Then given a similarity $\varphi\in\Phi_\infty(\cA)$ there exists
a normalized isomorphism $f\in\iso(\cA,\varphi)$ such that  $f^S\in\aut(S)$ for all $S\in\frS_0(\cA)$.
\ecrllr
\proof Let $\varphi\in\Phi_\infty(\cA)$.  Then $\varphi$ is induced by a normalized isomorphism $f$. 
Besides, by Theorem~\ref{270514a} without loss of generality we can
assume that there are no singular classes of $\cA$. Thus the required statement follows from the lemma below.

\lmml{200115x}
Let $\cA$ be a circulant S-ring without singular classes, and $f\in\iso(\cA)$. Then $f^S\le\aut(S)$ for 
all $S\in\frS_0(\cA)$.
\elmm
\proof Let $S\in\frS_0(\cA)$. By statement~(2) of Lemma~\ref{170614a}  we can assume that $S$ is principal, and 
hence $\rad(\cA_S)=1$.  By \cite[Theorem~4.1]{EP12} in this case $S=S_0\times\cdots\times S_k$ where 
$S_i\in\frS(\cA)$ for all~$i$, the S-ring $\cA_{S_0}$ is normal and  $\cA_{S_i}$ is an S-ring of rank $2$ and prime 
degree  for $i\ge 1$. Then $f^{S_0}\le\aut(S_0)$ by normality and so it suffices to verify that 
\qtnl{200115f}
f^{S_i}\le\aut(S_i),\quad i\ge 1.
\eqtn
However, since $\cA$ has no singular classes, from  \cite[Theorem~5.1]{EP03} it follows that any primitive $\cA$-section 
is projectively equivalent to a subnormal $\cA$-section. But for $i\ge 1$ the section $S_i$ has prime degree, and hence 
is primitive. Thus any such section is projectively equivalent to a subnormal $\cA$-section. Now, \eqref{200115f}
follows from statement~(2) of Lemma~\ref{170614a}.\bull\medskip

{\bf Proof of Theorem~\ref{071014a}.} The "only if\," part is obvious. To prove the "if\," part suppose that every 
quasidense S-ring over $G$ is separable. Let $\cA$ be an  arbitrary S-ring over $G$. First, we observe that
given a section $S\in\frS(\cA)$ belonging to a singular class, the $S$-extension of $\cA$ is separable if and only if
so is $\cA$ (Corollary~\ref{280514f}). However,  by statement~(1) of Lemma ~\ref{101213a} the rank of this
extension is less than the rank of~$\cA$. So without loss of
generality we can assume that the S-ring $\cA$ has no singular classes, or equivalently that $\cA$ is a quasidense
S-ring. Thus $\cA$ is separable by the assumption.\bull

\section{Proof of Theorems~\ref{280414a} and \ref{170113a}}

Before proving the theorems we give the exact definition of $\cA$-multipliers discussed in Introduction,
and study them in some details.
Let $\cA$ be an S-ring over a cyclic group $G$. Suppose that for every section $S\in\frS_0(\cA)$ we are given an
automorphism $\sigma_S\in\aut(S)$.

\dfntnl{200313b}
The family $\fS=\{\sigma_S\}_{S\in\frS_0(\cA)}$ is called an  $\cA$-multiplier if the following two 
conditions are satisfied for all sections $S_1,S_2\in\frS_0(\cA)$:
\nmrt
\tm{SM1} if $S_1\succeq S_2$, then $(\sigma_{S_1})^{S_2}=\sigma_{S_2}$,
\tm{SM2} if $S_1\sim S_2$, then $(\sigma_{S_1})^{f_{S_1,S_2}}=\sigma_{S_2}$
where $f_{S_1,S_2}$ is the projective isomorphism from $S_1$ onto $S_2$.\footnote{The equality statement here
is equivalent to saying that $\sigma_{S_1}$ and $\sigma_{S_2}$ are induced by raising to the same power.}
\enmrt
\edfntn

Clearly, the  set of all  multipliers of $\cA$  forms a group; we denote it by $\mult(\cA)$.
The restriction $\fS^T$ of a  multiplier $\fS\in\mult(\cA)$ to a section $T\in\frS(\cA)$,
is defined in the same way as for outer multipliers. \medskip

When $\cA$ is a coset S-ring, the sets $\mult(\cA)$ and $\fmult(\cA)$ are closely
related. Namely, in this case $\cA_S=\mZ S$  for all $S\in\frS_0(\cA)$ (Theorem~A8.3), and hence 
$\aut_\cA(S)=1$ for such~$S$. 
Therefore if $\fS=\{C_S\}$ is an outer $\cA$-multiplier, then $C_S$ is a singleton consisting of
an automorphism $\sigma_S\in\aut(S)$  for all~$S$,
and hence  $\fS':=\{\sigma_S\}$ is an  $\cA$-multiplier.  Clearly, the mapping $\fS\mapsto\fS'$
is a group isomorphism from $\fmult(\cA)$ onto $\mult(\cA)$.\medskip

Let now $\cA$ be a quasidense S-ring. Then $\frS_0(\cA_0)\supset\frS_0(\cA)$ where $\cA_0$ is
the coset closure of~$\cA$. Since $\cA_0$ is a coset S-ring, the discussion in the previous
paragraph shows that every $\cA_0$-multiplier $\fS_0$ induces an ordinary $\cA_0$-multiplier $\fS'$. The
restriction of it to $\frS_0(\cA)$ produces a family $\fS$ which is obiously an outer $\cA$-multiplier.
The mapping from $\fmult(\cA_0)$ to $\mult(\cA)$ that takes $\fS_0$ to $\fS$, is
denoted by $\eta$.

\lmml{260814c}
The mapping $\eta:\fmult(\cA_0)\to\mult(\cA)$ is a group isomorphism.
\elmm
\proof Let $S\in\frS_0$ where $\frS_0=\frS_0(\cA_0)$. Then $S_p\in\frS_0(\cA)$ for all prime $p$ dividing 
$|S|$ (Corollary A10.10) where $S_p$ is the Sylow $p$-subgroup of $S$. 
Since any element of $\aut(S)$ is uniquely determined by its projections to $\aut(S_p)$, an outer
$\cA_0$-multiplier $\fS_0$ is uniquely determined by the $\cA$-multiplier $\eta(\fS_0)$. This shows
that the mapping $\eta$ is injective. To prove that it is surjective, let 
$$
\fS=\{\sigma'_S\}_{S\in\frS_0(\cA)}
$$ 
be an $\cA$-multiplier. Given $S\in\frS_0$ set   $\sigma_S:=\prod_p\sigma'_{S_p}$ 
 (here $S_p\in\frS_0(\cA)$). Then $\sigma_S\in\aut(S)$ for all $S\in\frS_0$, and
$\sigma^{}_S=\sigma'_S$ for $S\in\frS_0(\cA)$. Let us check that conditions (M1) and (M2) are
satisfied for the family
$$
\fS_0=\{C_S\}_{S\in\frS_0}.
$$
 where $C_S=\{\sigma_S\}$.\medskip

Suppose first that $S,T\in\frS_0$ and $T\preceq S$.
Then $T_p\preceq S_p$ for all $p$.  By condition (SM1) for $\fS$ this implies that
$$
\sigma_T=\prod_p\sigma'_{T_p}=\prod_p (\sigma'_{S_p})^{T_p}=(\sigma_S)^T,
$$
as required. To verify (M2) let $S,T\in\frS_0$ and $S\sim T$. Then $S_p\sim T_p$ for all $p$. Therefore
by Corollary~\ref{270814a} we have $T_p=(S_p)^{f_{S_p,T_p}}$ and the restriction of $f_{S,T}$ to $S_p$
equals $f_{S_p,T_p}$. By condition (SM2) for $\fS$ this implies that
$$
\sigma_T=\prod_p\sigma'_{T_p}=\prod_p (\sigma'_{S_p})^{f_{S_p,T_p}}=(\sigma_S)^{f_{S,T}}.
$$
Thus $\fS_0$ is an outer $\cA_0$-multiplier. Since obviously $\eta(\fS_0)=\fS$, we are done. 
\bull\medskip

{\bf Proof of Theorem~\ref{280414a}.}
The S-ring $\cA_0$ is a coset one, and so separable (Theorem~A9.1).  Therefore 
$\Phi(\cA_0)^\cA$ is a subgroup of the group $\Phi_\infty(\cA)$.
Conversely, let $\varphi\in \Phi_\infty(\cA)$. Then by Corollary~\ref{280514a} there exists a normalized isomorphism
$f\in\iso(\cA,\varphi)$ such that   $f^S\in\aut(S)$ for all $S\in\frS_0(\cA)$. Set
$\fS=\{\sigma_S\}_{S\in\frS_0(\cA)}$
where $\sigma_S=f^S$. Then condition (SM1) is trivially satisfied whereas condition (SM2) follows from
statement~(2) of Lemma~\ref{170614a}. Thus $\fS$ is an $\cA$-multiplier.  By Lemma~\ref{260814c}  it can be extended  to an outer $\cA_0$-multiplier $\fS_0$.   Denote by $\varphi_0$ the similarity of the
S-ring $\cA_0$ that corresponds to $\fS_0$ (Theorem~\ref{071212d}). Since $\fS_0$ extends $\fS$,
the restriction of $(\varphi_0)_S$ to $\cA_S$ equals  the similarity $\varphi_S$ for all $S\in\frS_0(\cA)$. Therefore 
the similarities $(\varphi_0) ^\cA$ and $\varphi$ coincide on  all principal $\cA$-sections. By Lemma A3.1 this implies 
that $(\varphi_0) ^\cA=\varphi$ as required.\bull\medskip

{\bf Proof of Theorem~\ref{170113a}.} To prove the "only if\," part, suppose that the S-ring $\cA$ is separable.
Let $\fS'\in\fmult(\cA)$. Then by Theorem~\ref{071212d} there exists $\varphi\in\Phi(\cA)$ such that
$\fS'=\fS(\varphi)$. So by Theorem~\ref{280414a} there exists $\varphi_0\in\Phi(\cA_0)$ such that
$\varphi=(\varphi_0)^\cA$. Again by Theorem~\ref{071212d} applied this time to $\cA_0$, we find an outer
$\cA_0$-multiplier $\fS_0=\{C_S\}_{S\in\frS_0}$ corresponding to $\varphi_0$. Since $\cA_0$ is
a coset S-ring, the coset $C_S$ is a singleton consisting of an automorphism $\sigma_S\in\aut(S)$.
Therefore  $\fS'=\theta(\fS)$ where $\fS=\{\sigma_S\}_{S\in\frS_0(\cA)}$. Since obviously
$\fS\in\mult(\cA)$, we are done.\medskip

To prove the "if\," part, suppose that the homomorphism $\theta$ is surjective. To check that
the S-ring $\cA$ is separable, let $\varphi\in\Phi(\cA)$. Then by the hypothesis  there exists
an $\cA$-multiplier $\fS$ such that $\theta(\fS)=\fS(\varphi)$.  Set 
$$
\fS_0=\eta^{-1}(\fS)
$$
where $\eta$ is the isomorphism in Lemma~\ref{260814c}.  By Theorem~\ref{071212d}
there exists $\varphi_0\in\Phi(\cA_0)$ such that $\fS_0=\fS(\varphi_0)$. Now, since the S-ring
$\cA_0$ is separable (Theorem~A9.1), the similarity $\varphi_0$ is induced by an isomorphism $f\in\iso(\cA_0)$.
Since obviously $\varphi=(\varphi_0)^\cA$, we have $f\in\iso(\cA,\varphi)$, as required.\bull

\section{Proof of Theorem~\ref{080914a}}

By Theorem~\ref{071014a} we can restrict ourselves to quasidence S-rings. To prove that every quasidence
S-ring $\cA$ over a cyclic $p$-group $G$ is separable, we use induction on $|G|$.  If $\rad(\cA)=1$, then
$\cA$ is the tensor product of S-rings of rank~$2$ and a normal S-ring \cite[Theorem~4.1]{EP12}. 
Since $G$ is a cyclic $p$-group, this implies that $\rk(\cA)=2$ or $\cA$ is normal. Thus  $\cA$ is  separable: 
in the former case this is obvious whereas in the latter one this follows from \cite[Theorem~6.6]{EP01ce}. 
Let now $\rad(\cA)\ne 1$.

\lmml{071014b}
Let $\cA$ be an S-ring over a cyclic $p$-group. Suppose that $\rad(\cA)\ne 1$.  Then
$\cA$ is a proper $S$-wreath  product  such that $\rad(\cA_S)=1$
 or $|S|=4$.
\elmm
\proof By Theorem~A5.2 $\cA$ is a proper $S$-wreath  product where $S=U/L$. Without loss of generality we
assume  that the section $S$ is of minimal possible order. Then there exists $X\in\cS(\cA)_{G\setminus U}$
such that $\rad(X)=L$ (here we use that the set of $\cA$-groups is linearly ordered by inclusion). Therefore the $\cA$-section 
$T=\grp{X}/\rad(X)$ is of trivial radical and contains $S$ as a subsection.  This implies that 
$\rk(\cA_T)=2$ or $\cA_T$ is normal (see above). In the former case, obviously $\rad(\cA_S)=1$. In the
latter case, $\cA_T$ is a cyclotomic S-ring with trivial radical  \cite[Lemma~7.1]{EP01ce}.
However, from \cite[Theorem~7.3]{EP11} it  follows  that in a cyclotomic S-ring with trivial radical over a cyclic 
$p$-group, any of its subsections of order $\ne 4$ is of trivial radical. Thus $\rad(\cA_S)=1$
 or $|S|=4$. as required.\bull\medskip

By Lemma~\ref{071014b} the S-ring~$\cA$ is a proper $S$-wreath
product where the section $S$ is normal or of trivial radical. By Theorem~\ref{170113a} it suffices to verify that
given an outer $\cA$-multiplier $\fS$ there exists an $\cA$-multiplier $\fS'$ such that $\theta(\fS')=\fS$. 
To do this let $S=U/L$. Denote by $\varphi$ the similarity of~$\cA$ corresponding to~$\fS$ (Theorem~\ref{071212d}).
Then by Lemma~\ref{281113d} we have
\qtnl{081014a}
\fS(\varphi_1)=\fS^U\qaq \fS(\varphi_0)=\fS^{G/L}
\eqtn
where $\varphi_1=\varphi_U$ and $\varphi_0=\varphi_{G/L}$. 
Since $\cA_U$ and $\cA_{G/L}$ are separable (the inductive hypothesis),
these similarities belong to $\Phi_\infty(\cA_U)$ and $\Phi_\infty(\cA_{G/L})$ respectively.  So by 
Corollary~\ref{280514a} there exist normalized isomorphisms $f_1\in\iso(\cA_U,\varphi_1)$ 
and $f_0\in\iso(\cA_{G/L},\varphi_0)$ such that   $(f_1)^T\in\aut(T)$ for all $T\in\frS_0(\cA_U)$ and
$(f_0)^T\in\aut(T)$ for all $T\in\frS_0(\cA_{G/L})$.  Then the families
\qtnl{031214a}
\fS_1=\{(f_1)^T\}_{T\in \frS_0(\cA_U)}\qaq
\fS_0=\{(f_0)^T\}_{T\in \frS_0(\cA_{G/L})}
\eqtn
are $\cA_U$- and $\cA_{G/L}$-multipliers respectively (see the proof of Theorem~\ref{280414a}). Moreover, 
by~\eqref{081014a}, they extend the multipliers $\fS^U$  and $\fS^{G/L}$.\medskip

Without loss of generality we can assume that $(\fS_1)^S=(\fS_0)^S$, or equivalently 
\qtnl{021214c}
(f_1)^S=(f_0)^S.
\eqtn
To see this, we first observe that  the permutation $\sigma:=(f_1^S)^{-1}f_0^S$ induce the trivial similarity 
of the S-ring $\cA_S$, because both $(f_1)^S$ and $(f_0)^S$ induce  $\varphi_S$. 
So $\sigma$ belongs to $\aut(\cA_S)$, and by the choice of~$S$ also to $\aut_\cA(S)$. Then there 
exists $f\in\aut(\cA_U)$ such that 
\qtnl{021214a}
f^S=\sigma\qaq f^T\in\aut_\cA(T)\  \text{for all}\ T\in\frS_0(\cA_U).
\eqtn
Indeed, the S-ring $\cA_U$ is quasidense and schurian.\footnote{The fact 
that any circulant S-ring over a $p$-group, is schurian, was proved by R.~P\"oschel for odd $p$, and by I.~Kov\'acs 
for $p=2$, see the discussion in~\cite{EP11}.} Therefore by \cite[Theorem~3.5 and Remark 3.6]{EP11}  
there exists $f\in\aut(\cA_U)$ such that $f^S=\sigma$ and $f^T\in\aut_\cA(T)$ for 
any $\cA$-section $T$  of  trivial radical. It follows that $f^T\in\aut_\cA(T)$ for all subprincipal
sections~$T$. Since each nontrivial class of projective equivalence in a cyclic $p$-group is a singleton, this is true also for all
quasisubprincipal sections. Thus, we found an automorphism~$f$ satisfying~\eqref{021214a}. 
Now, equality~\eqref{021214c} holds with $f_1$ replaced by $f_1f$.\medskip

To complete the proof, let $T\in\frS_0(\cA)$. Then $T$ belongs to $\frS_0(\cA_U)$ or $\frS_0(\cA_{G/L})$: 
indeed, this is true if $T$ is
a subprincipal section, and hence for any $T$ because each nontrivial class of projective equivalence in a cyclic $p$-group 
is a singleton. Moreover, by~\eqref{021214c} we have $(f_1)^T=(f_0)^T$
if  $T\preceq S$. Thus, the following element of the group $\aut(T)$ is well defined:
$$
\sigma'_T:=\css
(f_1)^T,       &\text{if $T\in \frS_0(\cA_U)$,}\\
(f_0)^T,       &\text{if $T\in\frS_0(\cA_{G/L})$.}\\
\ecss
$$
Then the family $\fS':=\{\sigma'_T\}_{T\in\frS_0(\cA)}$ is obviously an $\cA$-multiplier.
Furthermore, from~\eqref{081014a} and~\eqref{031214a}  it follows that the multipliers $\theta((\fS')^U)$ and 
$\theta((\fS')^{G/L})$ correspond to the similarities $(\varphi_1)_U$ and $(\varphi_0)_{G/L}$ respectively. 
Therefore $\theta(\fS')$ corresponds to $\varphi$, and hence $\theta(\fS')=\fS$, as required.\bull

\section{Proof of Theorem~\ref{170113b}}
Let $\cA$ be a circulant S-ring. By Theorem~\ref{071014a} we can assume that it is quasidense.
Below for $S\in\frS(G)$ and $C\subset\aut(S)$ we set  $\wh C=\{\wh\sigma:\ \sigma\in C\}$.
For the reader convenience we cite here Lemma~A13.3.

\lmml{280113a}
Let $\cA$ be a quasidense circulant S-ring. Then
\nmrt
\tm{1} $\frS_0(\wh\cA)=\{\wh S:\ S\in\frS_0(\cA)\}$,
\tm{2} $\wh{\cA_0}=\wh\cA_{\,0}$ and $\frS_0(\wh\cA_{\,0})=\{\wh S:\ S\in\frS_0(\cA_0)\}$,
\tm{3} $\aut_{\wh\cA}(\wh S)=\wh{\aut_\cA(S)}$ for all $S\in\frS(\cA)$.\bull
\enmrt
\elmm

For a family $\fS=\{C_S\}_{S\in I}$ of sets $C_S\subset\aut(S)$  where $I\subset\frS(G)$, let us define the
family $\wh\fS=\{C_{S'}\}_{S'\in\wh I}$  of sets $C_{S'}\subset\aut( S')$,
by $C_{S'}=\wh{C_S}$ with $\wh{S}=S'$  where $\wh I=\{\wh S:\ S\in I\}$.

\crllrl{290814a}
Let $\cA$ be a quasidense circulant S-ring and $\fS$ an outer $\cA$-multiplier. Then the
family $\wh\fS$ is an outer $\wh\cA$-multiplier.
\ecrllr
\proof Let $\fS=\{C_S\}_{S\in\frS_0(\cA)}$. Then from statement~(3) of Lemma~\ref{280113a} it follows that
$C_{\wh S}$ is a coset of the subgroup $\aut_{\wh\cA}(\wh S)$ in the group $\aut(\wh S)$ for all $S$. Moreover, 
by  Lemmas~\ref{010914e}
and~\ref{010914a}, conditions (M1) and (M2) are satisfied  for the family $\wh\fS$. Thus   $\wh\fS$ is an 
outer $\wh\cA$-multiplier.\bull\medskip

Let us come back to the proof of Theorem~\ref{170113b}. Let $\cA$ be a  separable quasidense circulant S-ring. By
duality we only have  to prove that the S-ring $\wh\cA$ is separable. However, the latter ring
is quasidense by statement~(2) of Theorem~A5.9.  Therefore by Theorem~\ref{170113a} and Lemma~\ref{260814c}
it suffices to verify that the homomorphism $\wh\eta\wh\theta$ is surjective, where $\wh\eta$ and $\wh\theta$
are defined  for the S-ring $\wh\cA$ as in Lemma~\ref{260814c} and in formula~\eqref{260814b}, respectively.\medskip

Let $\fS'=\{C_{S'}\}_{S'\in\frS_0(\wh\cA)}$ be an outer $\wh\cA$-multipler. Then
by Corollary~\ref{290814a} applied to $\wh\cA$ and $\fS'$,  the family $\fS:=\{C_S\}_{S\in\frS_0(\cA)}$
is an outer $\cA$-multiplier where $C_S$ is defined by $\wh{C_S}=C_{S'}$  with $S'=\wh{S}$
(here we identify the duals to  $\wh\cA$ and $\wh{ C_S}$  with $\cA$ and $C_S$ respectively).
Since the S-ring $\cA$ is separable, by Theorem~\ref{170113a} and Lemma~\ref{260814c}
there exists an outer $\cA_0$-multiplier $\fS_0$ such that $(\fS_0)^{\eta\theta}=\fS$. By Corollary~\ref{290814a}
applied to $\cA_0$ and $\fS_0$, the family $\wh{\fS_0}$ is an outer $\wh{\cA_0}$-multiplier, and
hence an outer $\wh\cA_{\,0}$-multiplier (statement~(2) of Lemma~\ref{280113a}). Thus the surjectiviity statement is true
because
$$
(\wh{\fS_0})^{\wh\eta\wh\theta}=
(\wh{(\fS_0)^{\eta}})^{\wh\theta}=
\wh{(\fS_0)^{\eta\theta}}=
\wh\fS=\fS'.
$$
Here the first equality follows from statement~(1) of Lemma~\ref{280113a} and the definition of $\wh\eta$,
whereas the third and fourth equalities follow from the definitions of $\fS_0$ and $\fS$ respectively. To prove the
second equality let $S'\in\frS_0(\wh\cA)$, and let $S\in\frS_0(\cA)$ be such that $\wh S=S'$.
Denote the $S'$-component of $\wh{(\fS_0)^{\eta}}$ by
$\sigma_{S'}$. Then by statement~(1) of Lemma~\ref{280113a} the $S'$-component of
$(\wh{(\fS_0)^{\eta}})^{\wh\theta}$ equals
$$
\aut_{\wh\cA}(S')\sigma_{S'}=\wh{\aut_\cA(S)}\wh{\sigma_S}
$$
which is the $S'$-component of $\wh{(\fS_0)^{\eta\theta}}$ as required.\bull

\end{document}